\documentclass[preprint,12pt]{elsarticle}

\usepackage{amsfonts, amsmath, amscd}
\usepackage[psamsfonts]{amssymb}

\usepackage{amssymb}

\newcommand{\pr}{\par {\bf Proof.} }   
\newcounter{df}
\newenvironment{df}{\par
\refstepcounter{df}
{\bf Definition \arabic{df}.} }{}
\newcounter{pro}
\newenvironment{pro}{\par
\refstepcounter{pro}
{\bf Proposition \arabic{pro}.} }{}
\newcounter{rem}
\newenvironment{rem}{\par
\refstepcounter{rem}
{\bf Remark \arabic{rem}.} }{}
\newcounter{exa}
\newcounter{teo}
\newenvironment{teo}{\par
\refstepcounter{teo}
{\bf Theorem \arabic{teo}.} \it }{}
\newcounter{cor}
\newcounter{st}
\newcounter{lem}
\newenvironment{lem}{\par
\refstepcounter{lem}
{\bf Lemma \arabic{lem}.} \it }{}
\makeatletter                                                   %
\renewcommand{\section}{\@startsection{section}{1}
{\parindent}{3.5ex plus 1ex minus 0.2ex}{2.3ex plus 0.2ex}{\bf}}
\makeatother                                                    %

\begin{document}

\begin{frontmatter}

\title{Groups of quasi-invariance and the Pontryagin duality\tnoteref{label1}}\tnotetext[label1]{The author was partially supported  by Israel Ministry of Immigrant Absorption}

\author{S.S. Gabriyelyan}
\ead{saak@math.bgu.ac.il}

\address{Department of Mathematics, Ben-Gurion University of the
Negev, Beer-Sheva, P.O. 653, Israel}

\begin{abstract}
A Polish group $G$ is called a group of quasi-invariance  or a QI-group, if there
exist a locally compact group $X$ and a probability measure $\mu$
on $X$ such that 1) there exists a continuous monomorphism of $G$ to $X$, and 2) for each $g\in X$ either $g\in G$ and the shift $\mu_g$ is equivalent to $\mu$ or $g\not\in G$ and $\mu_g$ is orthogonal to $\mu$. It is proved that $G$ is a $\sigma$-compact subset of $X$.  We show that there exists a quotient group $\mathbb{T}^H_2$ of $\ell^2$ modulo a discrete subgroup which is a Polish monothetic non locally quasi-convex (and hence nonreflexive) pathwise connected QI-group, and such that the bidual of $\mathbb{T}^H_2$ is not a QI-group. It is proved also that the bidual group of a QI-group may be not a saturated subgroup of $X$.
\end{abstract}


\end{frontmatter}


\section{Introduction}
\label{}

Let $X$ be a Polish group and ${\cal B}$ the family of its Borel sets. Let $E\in {\cal B}$. The image and the inverse image of $E$ are denoted by $g\cdot E$ and $g^{-1} E$ respectively. Let $\mu$ and $\nu$ be probability measures on $X$.  We write $\mu\ll\nu $ ($\mu\sim\nu , \mu\perp\nu$) if $\mu$ is absolutely continuous relative to $\nu$ (respectively: equivalent, mutually singular).  For $g\in X$ we denote by $\mu_g$ the measure
determined by the relation $\mu_g (E) := \mu (g^{-1} E)$, $\forall E$. The set of all $g$ such that $\mu_g \sim\mu$ is denoted by $E(\mu)$. The Mackey-Weil Theorem asserts that $X=E(\mu
)$ for some $\mu$ iff $X$ is locally compact. Some algebraic and topological properties of $E(\mu)$ are considered in \cite{Ga1} and \cite{Ga2}. In particular, it was proved that  $E(\mu)$ always admits a Polish group topology and, as a subgroup of $X$, is a $G_{\delta\sigma\delta}$-set.
The Polish group topology is defined by the strong operator topology  in the following way. If $g,h\in E(\mu)$, and $\left\{\mu^n \right\}$ is a countable dense subset in $L^1 (\mu
)$ (with $\mu^1 = \mu$), then the following metric on $E(\mu )$
\[
d(h,g) = \sum_{n=1}^{\infty} \frac{1}{2^n } \left( \frac{\|
\mu^n_g - \mu^n_h  \| }{1+\| \mu^n_g - \mu^n_h \| } + \frac{\| \mu^n_{g^{-1}} - \mu^n_{h^{-1}} \| }{1+\| \mu^n_{g^{-1}} - \mu^n_{h^{-1}} \| }  \right) .
\]
defines the Polish group topology which is finer than the topology induced from $X$. Note that although the metric $d$ depends on the chosen sequence $\left\{\mu^n \right\}$, the Polish group topology is unique and does not depend on $d$.

For a topological group $G$, the group $G^{\wedge}$ of continuous homomorphisms (characters) into the torus $\mathbb{T} =\{ z\in \mathbb{C} : \; |z|=1\}$ endowed with the compact-open topology is called  the {\it character group} of $G$ and $G$ is named {\it Pontryagin reflexive} or {\it reflexive} if the canonical homomorphism $\alpha_G : G\to G^{\wedge\wedge} , g\mapsto (\chi\mapsto (\chi, g))$ is a topological isomorphism. A subset $A$ of $G$ is called {\it quasi-convex} if for every $g\in G\setminus A$, there is some $\chi \in A^{\triangleright} := \{ \chi \in G^{\wedge} : \; {\rm Re} (\chi , h) \geqslant 0, \forall h\in A\}$, such that ${\rm Re} (\chi, g) <0$, \cite{Vil}. An Abelian topological group $G$ is called {\it locally quasi-convex} if it has a neighborhood  basis of the neutral element $e_G$, given by quasi-convex sets. The dual $G^{\wedge}$ of any topological Abelian group $G$ is locally quasi-convex \cite{Vil}. In fact, the sets $K^{\triangleright}$, where  $K$ runs through the compact subsets of $G$, constitute a neighborhood basis of $e_{G^{\wedge}}$ for the compact open topology.

Following E.G.Zelenyuk and I.V.Protasov \cite{ZP1}, we say that a sequence $\mathbf{u} =\{ u_n \}$ in a group $G$ is a $T$-{\it sequence} if there is a Hausdorff group topology on $G$ for which $u_n $ converges to zero. The group $G$ equipped with the finest group topology with this property is denoted by $(G, \mathbf{u})$.

 Set
\[
\mathbb{Z}^{\infty}_b = \left\{ \mathbf{n} =(n_1,\dots, n_k, n_{k+1} ,\dots)  | \; n_i \in \mathbb{Z} \mbox{ and } |\mathbf{n} |_b := \sup \{ | n_i |\} <\infty \right\},
\]
\[
\mathbb{Z}^{\infty}_0 =\{ \textbf{n}= (n_1,\dots, n_k, 0,\dots) | n_j \in \mathbb{Z} \}.
\]
We will consider the spaces $l^p$ and $c_0$. For our convenience we set $l^0 :=c_0$. Evidently, $\mathbb{Z}_0^{\infty}$ is a closed discrete subgroup of $l^p$ for any $0\leq p<\infty$.

The following groups play a crucial role in our consideration
\[
\mathbb{T}^H_p := \left\{ \omega =(z_n )\in \mathbb{T}^{\infty} \, |
\quad \sum_{n=1}^{\infty} |1-z_n |^p <\infty \right\} , 0<p<\infty,
\]
\[
 \mathbb{T}^H_0 := \left\{ \omega =(z_n )\in \mathbb{T}^{\infty} \, |  \quad z_n \to 1 \right\} .
\]
It is easy to prove that $\mathbb{T}^H_p $ are Polish groups with pointwise multiplication and the topology generated by the metric
\[
d_p (\omega_1 , \omega_2 )= \left( \sum_{n=1}^{\infty} |z_n^1
-z_n^2 |^p \right)^{\min (1,\frac{1}{p}) } , \mbox{ if } 0<p<\infty , \; \mbox{ and }
\]
\[
 d_0 (\omega_1 , \omega_2 )= \sup ( |z_n^1 -z_n^2 |, n=1,2,\dots ),  \mbox{ if } p=0.
\]
We also need more complicated groups, which are defined in \cite{AaN}.
Suppose that $a_n \geq 2$ are integers such that $\sum_{n=1}^{\infty} \frac{1}{a_n} <\infty$ and set $\gamma (1)=1, \gamma (n+1)=\Pi_{k=1}^n a_k (n\geq 1)$. If $1\leq p <\infty$ and $z\in \mathbb{T}$, we set
\[
\| z \|_p = \left( \sum_{n=1}^{\infty} | 1 - z^{\gamma (n)} |^p \right)^{1/p} \mbox{ and }.
\]

Put $G_p = \{ z\in \mathbb{T} : \, \| z\|_p <\infty \}$ and $Q= \{ z\in \mathbb{T} ; \, z^{\gamma (n)} =1$ for some $n\}$. Then $(G_p , \| .\|_p )$ is a Polish group and $Q$ is its dense subgroup \cite{AaN}. It is clear that all $G_p$ are totally disconnected.

For $p=0$ and $z\in \mathbb{T}$ we set $G_0 = \{ z\in \mathbb{T} : \,  z^{\gamma (n)} \to 1 \}$ and $ \| z \|_0 = \sup \{ |1 -z^{\gamma (k)} |, k=1,2\dots \}$. Then $(G_0 , \| .\|_0 )$ is a Polish group \cite{Ga3}.

This article was inspired by the following old and important problem in abstract harmonic analysis and topological algebra: find the ``right'' generalization of the class of locally compact groups. The commutative harmonic analysis gives us one of the best indicator for the ``right'' generalization - the Pontryagin duality theorem. The Pontryagin theorem is known to be true for several classes of non locally compact groups: the additive group of a Banach space, products of locally compact groups, complete metrizable nuclear groups \cite{Ban}, \cite{Kap} \cite{Smi}. These examples suggest searching for possible generalizations not only in the direction of duality theory. The existence of the Haar measure plays a crucial role in harmonic analysis. As it was mentioned above, existence of a left (quasi)invariant measure is equivalent to the local compactness of a group. Therefore we can use some similar notion only. Groups of the form $E(\mu)$ are natural candidates. On the other hand, if a probability measure $\mu$ on $\mathbb{T}$ is ergodic under $E(\mu)$, then for each $g\in \mathbb{T}$ either $\mu_g \sim \mu$ or $\mu_g \perp \mu$ \cite{HMP}. Therefore, taking into consideration ergodic decomposition, we propose the following generalization.

\begin{df} \label{d1}
{\it A Polish group $G$ is called a group of quasi-invariance  or a QI-group, if there
exist a local compact group $X$ and a probability measure $\mu$
on $X$ such that $G$ is continuously embedded in $X$, $E(\mu)=G$ and $\mu_g \perp
\mu$ for all $g\not\in E(\mu)$.}
\end{df}

We will say that  $G$ is represented in $X$ by $\mu$ and denote it by $E_{\mu}$. It is clear that, if $G$ is Abelian, then we can assume that $X$ is compact.
In the general case we can define QI-groups in the following way (taking into account Theorems 5.14 and 8.7 in \cite{HR}): a topological group $G$ is called  a QI-group if there exists a compact normal subgroup $Y$ such that $G/Y$ is a Polish QI-group. In the article we are restricted to the Abelian Polish case only. By the definition, it is clear that a QI-group has enough continuous characters.

Evidently, any locally compact Polish group is a QI-group. Let ${\cal H}$ be a separable real Hilbert space. Then ${\cal H}$ is a QI-group, since it is $E(\mu )$ for a Gaussian measure on $\mathbb{R}^{\infty}$ \cite{ShF} and we can continuously embed  $\mathbb{R}^{\infty}$ into $\left( \mathbb{T}^2 \right)^{\infty} = \mathbb{T}^{\infty}$ in the usual way. Moreover, any $l^p, 0<p\leq 2,$ is a QI-group \cite{ChM}. In \cite{AaN}, the authors proved that $G_1$ and $G_2$ are QI-groups (see also \cite{Par}). We will give a simple straightforward proof that $\mathbb{T}^H_1$ and $\mathbb{T}^H_2$ are QI-groups too (cf. \cite{Sad} for $\mathbb{T}^H_2$, see also \cite{Hor}). If a probability measure $\mu$ on $\mathbb{T}$ is ergodic under $E(\mu)$, then $E(\mu) =E_{\mu}$ is a QI-group \cite{HMP}. Hence, in the category of Polish groups, the set ${\cal GQI}$ of all groups of quasi-invariance is wider than the class of locally compact groups.

Choosing of such groups is motivated not only by the above-mentioned. Let $\mu$ on $\mathbb{T}$ be ergodic under $E(\mu) =E_{\mu}$. J.Aaronson and M.Nadkarni \cite{AaN} showed that $E_{\mu}$ is the eigenvalue group of some non-singular transformation and illustrated a basic interaction between eigenvalue groups and $L^2$ spectra. Moreover, they computed the Hausdorff dimension of some $E_{\mu}$, which is important in connection with the dissipative properties of a non-singular transformation \cite{Aar}. A deep property of the eigenvalues of the action of $E_{\mu}$ gives us a key property for solving subtle problems about spectra of measures around the Wiener-Pitt phenomenon. These and other applications to harmonic analysis are given in \cite{HMP} and \cite{PAR}.
Below we prove that a QI-group is even a $\sigma$-compact subgroup of some locally compact group. Hence, on the one hand, QI-groups play a very important role in non-singular dynamics and harmonic analysis and, on the other hand, they are not ``very big'' (as, for example, the unitary group $U(H)$ of the separable Hilbert space or the infinite symmetric groups $S_{\infty}$).  These arguments allow us to consider the notion ``to be a QI-group'' as a possible generalization of  the notion ``to be a locally compact group'' and explain our interest in such groups.

The main goal of the article is to consider some general problems of the Pontryagin duality theory for groups of quasi-invariance. The following question is natural:
\begin{enumerate}
\item[] {\bf Question 1.} {\it Are all groups of quasi-invariance Pontryagin reflexive}?
\end{enumerate}
It is clear that locally compact Polish groups and $l^p, 1\leq p\leq 2,$ are reflexive. We prove that $\mathbb{T}^H_1$ is reflexive too.  On the other hand, any group $l^p, 0< p<1,$ is not even locally quasi-convex (8.27,\cite{Aus}) and, hence, not reflexive. Since the bidual group of a Polish group is always locally quasi-convex and Polish \cite{Cha}, we can ask the following.
\begin{enumerate}
\item[] {\bf Question 2.} {\it Is the bidual $G^{\wedge\wedge}$ of a QI-group $G$ a QI-group}?
\end{enumerate}
The answer on question 2 is also negative. We prove that the bidual group of $\mathbb{T}^H_2$ is not a QI-group.  We do not know the answer on the next question.
\begin{enumerate}
\item[] {\bf Question 3.} {\it Let $G$ be a locally quasi-convex {\rm QI}-group. Is $G$ reflexive}?
\end{enumerate}

In \cite{HMP}, the authors proved that each QI-group is saturated.
Since the bidual of a QI-group $G$ is not always a QI-group, we can ask whether the bidual group is saturated. It is proved that  $G_2^{\wedge\wedge} = G_0$ (and, hence,  $G_2$ is not reflexive). Since $G_0$ is not a saturated subgroup of $\mathbb{T}$ \cite{HMP}, this shows that the bidual group of a QI-group may be not saturated.

On the other hand, the groups $\mathbb{T}^{H}_p$ are interesting from the general point of view of the Pontryagin duality theory. We  prove the following.
\begin{itemize}
\item If $1< p <\infty$, then $\mathbb{T}^H_p$ is a {\it monothetic non} locally quasi-convex {\it Polish} group and, hence, non reflexive. V.~Pestov \cite{Pes} asked whether every $\check{{\rm C}}$ech-complete group $G$ with sufficiently many characters is a reflexive group. Hence $\mathbb{T}^H_p$ gives another negative answer on this question (in 11.15 \cite{Aus} an even stronger counterexample is given).
\item If $p=0$ or $p=1$, then $\mathbb{T}^H_p$ is a {\it monothetic reflexive} Polish {\it non} locally compact group.
\item $\mathbb{T}^H_p$ is topologically isomorphic to $l^p /\mathbb{Z}^{\infty}_0$. Since $l^p$ is Pontryagin reflexive and $\mathbb{Z}^{\infty}_0$ is its closed discrete (and hence) locally compact subgroup, we see that their quotient  $\mathbb{T}^H_p$ is not locally quasi-convex. Thus the answer on question 14 \cite{Bru} is negative.
\item $\left( \mathbb{T}^H_p \right)^{\wedge\wedge}$ is topologically isomorphic to $c_0 / \mathbb{Z}^{\infty}_0$ and reflexive. Hence neither ``to be dual'' nor ``to be Pontryagin reflexive'' is not a three space property.
\item If $1<p<\infty$, then $\left( \mathbb{T}^H_p\right)^{\wedge} = \left( \mathbb{T}^H_0\right)^{\wedge} =\mathbb{Z}_0^{\infty}$ and is {\it reflexive}. Thus there exists a continual chain (under inclusion) of {\it Polish monothetic non} locally quasi-convex pathwise connected groups with the same {\it countable reflexive} dual. In fact, $\left( \mathbb{T}^H_0\right)^{\wedge}$ is a Graev free topological abelian group over the convergent sequence.
\end{itemize}

Analogous properties hold for the family of groups $G_p$. Let us note only that, actually, we can consider the group $G_p$ as a dually closed and dually embedded (totally disconnected) subgroup of $\mathbb{T}^H_p$.

\section{Main Results}

As it was mentioned above, $E(\mu)$ is a $G_{\delta\sigma\delta}$-subset of $X$.  For a QI-group we can prove the following.

\begin{pro} \label{p1}
{\it If a {\rm QI}-group $G$ is represented in $X$, then it is $\sigma$-compact in $X$.}
\end{pro}
\pr Since the Polish group topology $\tau$ on $G$ is unique, we can consider $G$ as $(E(\mu), d)$ for some probability measure $\mu$ on $X$. Since $\tau$ is finer than the topology on $X$, we can choose $\varepsilon_0 >0$ such that the $\varepsilon_0$-neighborhood $U_{\varepsilon_0}$ of the unit $e$ is contained in a compact neighborhood of $e$ in $X$. Thus ${\rm Cl}_X
U_{\varepsilon}$ is compact in $X$ for any $\varepsilon<\varepsilon_0$. If $\{ h_n\}$ is a dense countable subset of $(E(\mu), d)$, then $E(\mu)=\cup_n h_n U_{\varepsilon}$ for every $\varepsilon>0$. Therefore, if we will prove that ${\rm Cl}_X U_{\varepsilon} \subset E(\mu)$ for an enough small $\varepsilon<\varepsilon_0$, then $E(\mu)=\cup_n h_n {\rm Cl}_X U_{\varepsilon}$ is a $\sigma$-compact subset of $X$. Set $\varepsilon =\min (\varepsilon_0, 0.1)$. We will prove that ${\rm Cl}_X
U_{\varepsilon} \subset E(\mu)$. Let $g_n \to t$ in $X, g_n \in
U_{\varepsilon}.$ Assume the converse and $t\not\in E(\mu )$, i.e.
$\mu_t \perp\mu$. Choose a compact $K$ such that $\mu (K) > 0,9$ and $\mu_t (K)=0$. Choose a neighborhood $V_{\delta} (K)$ of $K$
such that $\mu_t (V_{\delta} (K)) < 0,1$.
Then there exists an integer $N$ such that
\begin{equation} \label{e1}
g_n^{-1} \cdot K \subset t^{-1} V_{\delta} (K), \forall n> N,
\mbox{ and } \mu_{g_n} (K)  = \mu (g_n^{-1}
K)<0,1 , \forall n> N.
\end{equation}
On the other hand, since $d(g,e) <\varepsilon $, then, by the definition of $d$, $\| \mu_g -
\mu\|<0,2$. But for any $g\in U_{\varepsilon}$ we have
\[
0,2 > \| \mu_g - \mu\| \geq \| \mu_g |_K - \mu |_K \|
\geq | \mu (g^{-1} K) - \mu (K) | ,
\]
\[
\mbox{and } \mu (g^{-1} K)= \mu(K) + (\mu (g^{-1} K) -\mu (K))
\geq 0,9- 0,2=0,7.
\]
In particular, $\mu_{g_n} (K) > 0,7$. This inequality  contradicts to (\ref{e1}). $\Box$

\begin{pro} \label{pro2}
{\it $\mathbb{T}^H_1$ and $\mathbb{T}^H_2$ are {\rm QI}-groups.}
\end{pro}
\pr We can capture the idea of how to construct of examples as follows. Let $\mu$ be absolutely continuous relative to the Haar measure $m_{\mathbb{R}}$ and assume that its density $f(x)$ is smooth. Then
\[
P(\varphi):= \int \sqrt{f(x)f(x+\varphi)} dx= \int f(x) \sqrt{1+\frac{1}{f(x)} (f(x+\varphi)- f(x))} dx =
\]
\[
1+\frac{\varphi}{2} \int f'(x) dx -\frac{\varphi^2}{8} \int \frac{ (f')^2 - 2ff''}{f} dx +
\]
\[
\frac{\varphi^3}{48} \int \frac{f^2 \cdot f''' - 6f f' f'' +3 (f')^3}{f^2} dx + O(\varphi^4).
\]
Hence we can expect: if $f$ is linear, then $P(\varphi) \sim 1+c\varphi$; and if $\int f' dx =0$, then  $P(\varphi) \sim 1-c\varphi^2$.

We identify $\mathbb{T}$ with $[-\frac{1}{2}; \frac{1}{2}), t\mapsto
\textrm{e}^{2i\pi t}$. Since $\alpha /2 < \sin \alpha <\alpha , \alpha\in (0, \pi /2)$, then for $\varphi\in [-\frac{1}{2} ; \frac{1}{2} ), z= {\rm e}^{2\pi i \varphi}, 0<p<\infty$, we have
\begin{equation} \label{2}
2^p \pi^p |\varphi |^p \geq | 1-z|^p = | 1- e^{2\pi i\varphi} |^p = 2^p |\sin \pi\varphi |^p \geq \pi^p |\varphi |^p , \; p>0.
\end{equation}

1) {\it Let us prove that $\mathbb{T}^H_1$ is a {\rm QI}-groups.}

Let $f(x) =x+1$. For
$\varphi\in [0;\frac{1}{2})$ we get
\[
f(x-\varphi ) =x+1- \varphi , \mbox{ if } x\in
[-\frac{1}{2} +\varphi;\frac{1}{2}), \mbox{ and }
\]
\[
f(x-\varphi ) =x+\frac{3}{2} -\varphi , \mbox{ if }  x\in [-\frac{1}{2}; -\frac{1}{2} +\varphi ].
\]
Then the routine computations give us the following
\[
P(\varphi ) \sim 1-\frac{8+5\sqrt{2}}{6+4\sqrt{2}} \varphi +O(\varphi^2)
\]
and $\max \{P(\varphi)\} = P(0)=1$ only at 0. Hence $P(\varphi )\to 1$ iff $\varphi\to 0$.
Consider the probability measures $\mu_n = f(x)
m_{\mathbb{T}} $ on $\mathbb{T}$. Set $\mu = \prod_n \mu_n$. Let $\omega =(z_n) = (\textrm{e}^{2i\pi \varphi_n })$. Then, by the Kakutani Theorem \cite{Ga1}, $\mu_{\omega} \not\perp\mu$ iff $\omega \in E(\mu)$ iff
\[
\prod_n P_n(\varphi_n ) <\infty \Leftrightarrow \sum_n \ln
P_{n} (\varphi_n ) <\infty \Leftrightarrow \sum_n | \varphi_n | <\infty .
\]
Since $\varphi_n \to 0$, then, by (\ref{2}), $|\varphi_n| \sim |1-
z_n | \cdot 2\pi$. Hence $\omega \in \mathbb{T}^H_1$.

2) {\it Let us prove that $\mathbb{T}^H_2$ is a {\rm QI}-groups.}

Let $f_c (x) =\frac{1}{a} \textrm{e}^{-
c|x|},$ where $a=\frac{2}{c} (1-\textrm{e}^{-c/2} )$. For
$\varphi\in [0;\frac{1}{2})$ we get
\[
f_c (x+\varphi ) =\frac{1}{a} \textrm{e}^{- c|x+\varphi |},  \mbox{ if } x\in
[-\frac{1}{2};\frac{1}{2}-\varphi ), \mbox{ and }
\]
\[
f_c (x+\varphi ) =\frac{1}{a}
\textrm{e}^{- c|x+\varphi -1|},  \mbox{ if } x\in [\frac{1}{2}-\varphi ;\frac{1}{2}).
\]
Then the simple computations give us
\[
P_c (\varphi ) = \int_{-\frac{1}{2}}^{\frac{1}{2}} \sqrt{f_c (x)f_c (x+\varphi
)} dx = \frac{1}{2} \textrm{sh}^{-1} \frac{c}{4} \left( 2\textrm{sh}
\frac{c}{4}(1-2\varphi ) + c\varphi \textrm{ch}\frac{c}{4}
(1-2\varphi ) \right).
\]
 It is easy to prove that
\begin{equation} \label{3}
1-\frac{1}{8} (c\varphi )^2 \leq P_c (\varphi ) \leq 1
-\frac{1}{32} (c\varphi )^2 , \quad \forall c\in [0;1] , \; \varphi \in [-\frac{1}{2};\frac{1}{2}).
\end{equation}
Consider the probability measures $\mu_n = f_{c_n} (x)
m_{\mathbb{T}} $ on $\mathbb{T}$. Set $\mu = \prod_n \mu_n$. Let $\omega =(z_n) = (\textrm{e}^{2i\pi \varphi_n })$. Then, by the Kakutani Theorem \cite{Ga1} and (\ref{3}) , $\mu_{\omega} \not\perp\mu$ iff $\omega \in E(\mu)$ iff
\[
\prod_n P_{c_n} (\varphi_n ) <\infty \Leftrightarrow \sum_n \ln
P_{c_n} (\varphi_n ) <\infty \Leftrightarrow \sum_n \left( c_n
\varphi_n \right)^2 <\infty .
\]
In particular, if $c_n =1$, then $\varphi_n \to 0$. Therefore, by (\ref{2}), $\varphi^2_n \sim |1-
z_n |^2 \cdot 4\pi^2$. Hence $\omega \in \mathbb{T}^H_2$. $\Box$

We do not know any general characterization of QI-vector spaces. In particular, we do not know the answer on the following question (taking into account Theorem 1 \cite{Hor})
\begin{itemize}
\item[] {\bf Question 4.} {\it Is there $p>2$ such that $l^p$ is a QI-group}?
\end{itemize}

The next proposition allows to prove the reflexivity of the dual group of a group $G$.

\begin{pro} \label{p3}
{\it Let $G$ be a Polish group. Set $H= {\rm Cl} (\alpha_G (G))$.
\begin{enumerate}
\item If $H=G^{\wedge\wedge} $, then $G^{\wedge}$ (and $G^{\wedge\wedge}$) is reflexive.
\item If $G$ is locally quasi-convex and $H=G^{\wedge\wedge} $, then $G$ is reflexive.
\end{enumerate} }
\end{pro}

\pr 1) Since $G$ is Polish and $G^{\wedge}$ is a $k$-space \cite{Cha}, then $\alpha_G$  and $\alpha_{G^{\wedge}}$ are continuous (corollary 5.12 \cite{Aus}). Since $H=G^{\wedge\wedge}$, $\alpha_{G^{\wedge}}$ is a continuous isomorphism. Let $\alpha_G^{\ast}$ be the dual homomorphism of $\alpha_{G}$. Then it is the converse to $\alpha_{G^{\wedge}}$, since
\[
\left( \alpha_G^{\ast} \circ \alpha_{G^{\wedge}} (\chi) , x\right) = \left( \alpha_{G^{\wedge}} (\chi), \alpha_G (x) \right) = \left( \chi , \alpha_G (x) \right) = (\chi, x), \; \forall \chi\in G^{\wedge}, \forall x\in G.
\]
Thus $\alpha_{G^{\wedge}}$ is a topological isomorphism and $G^{\wedge}$ is reflexive.

2) Let $G$ be locally quasi-convex. Then $\alpha_G (G)$ is an embedding with the closed image (Proposition 6.12 \cite{Aus}). Thus, if $H=G^{\wedge\wedge} $, then $G=G^{\wedge\wedge} $ is reflexive \cite{Cha}. $\Box$

\begin{pro} \label{p4}
{\it Let $0\leq p<\infty$. Then $\mathbb{T}^H_p$ is a monothetic Polish group which is topologically isomorphic to $l^p/ \mathbb{Z}_0^{\infty}$.}
\end{pro}

\pr 1) Let $\pi_p : l^p \to l^p/ \mathbb{Z}_0^{\infty}$
be the canonical map $0\leq p<\infty$. Denote by $\textbf{x}$ the class of equivalence of $(x_n)$, i.e. $\textbf{x} =(x_n )+\mathbb{Z}_0^{\infty}$. Then $\pi_p (x_n ) = \left( x_n (\textrm{mod}
1)\right)$. Indeed, $(x_n) +\mathbb{Z}_0^{\infty} = (y_n) +\mathbb{Z}_0^{\infty}$ iff there exists
an integer $N$ such that $y_n = x_n +m_n , m_n \in \mathbb{Z} ,
n=1,\dots ,N,$ and $y_n =x_n$ for $n>N$. Since $x_n$ and $y_n$
tend to zero, this is equivalent to $y_n =x_n (\textrm{mod}1)$. Set  $s_n = (y_n - x_n )
(\textrm{mod}1) \in [-\frac{1}{2} , \frac {1}{2} )$.
 Then the metric on $l^p/ \mathbb{Z}_0^{\infty}$ is defined as
\[
d^{\ast} (\textbf{x} , \textbf{y}) =\inf \left\{ d( (x'_n ), (y'_n
) ) , (x'_n )\in (x_n )+\mathbb{Z}_0^{\infty} , (y'_n )\in (y_n
)+\mathbb{Z}_0^{\infty} \right\},
\]
\begin{equation} \label{4}
\mbox{and so } d^{\ast} (\textbf{x} , \textbf{y}) = \left(\sum  |s_n |^p \right)^{\min(1,\frac{1}{p})}.
\end{equation}

Let $r: l^p / \mathbb{Z}_0^{\infty} \to \mathbb{T}^H_p , r
(\textbf{x}) = \left( \textrm{e}^{2\pi i \left( x_n (\textrm{mod}
1) \right)} \right)$. Evidently that $r$ is injective. If $(z_n) = (\textrm{e}^{2i\pi \varphi_n }) \in l^p , \varphi_n \in [-\frac{1}{2}; \frac{1}{2})$, then, by (\ref{2}), we have
\begin{equation} \label{5}
\pi^p \sum_{n=1}^{\infty} |\varphi_n |^p \leq \sum_{n=1}^{\infty} | 1-z_n |^p \leq 2^p \pi^p \sum_{n=1}^{\infty} |\varphi_n |^p .
\end{equation}
 Then
\begin{equation} \label{6}
d_p ( p(\textbf{x}), p(\textbf{y}) ) = \left( \sum
|\textrm{e}^{2\pi i s_n } -1 |^p \right)^{\min(1,\frac{1}{p})} =\left( \sum 2^p
|\sin \pi s_n |^p \right)^{\min(1,\frac{1}{p})} .
\end{equation}
Equations (\ref{4})-(\ref{6}) show that $\pi d^{\ast} \leq d
\leq 2\pi d^{\ast} $ and $r$ is surjective. Hence $r$ is a topological isomorphism.

Analogously, we can consider the case $p=0$.

2) S.~Rolewicz \cite{Rol} proved that $\mathbb{T}^H_0$ is monothetic. The case $0<p<\infty$ is considered analogically. We follow S.~Rolewicz \cite{Rol} and 1.5.4 \cite{DPS}.

Let $0<a_1 <1/2$ be an irrational number. For every $n>1$ we choose an irrational number $a_n$ such that
\begin{enumerate}
\item[a)] $a_1 >a_2>\dots >a_n>0$ are rationally independent.
\item[b)] $a_n < \frac{1}{2^n k_n}$, where $k_n$ is the smallest natural number such that for every $(n-1)$-tuple $(y_1,\dots, y_{n-1})$ of reals there exist integers $m_1,\dots, m_{n-1}$ and a natural number $k$ with $k\leq k_n$ and $|ka_s - y_s -m_s|<\frac{1}{2^n}$ for all $s=1,\dots n-1$ (the existence of such $k$ follows from the Kronecker Theorem).
\end{enumerate}
Now we set $\omega_0 =(z_n^0)$, where $z_n^0 ={\rm e}^{2\pi i a_n}$. Evidently $\omega_0 \in \mathbb{T}^H_p$ for every $0\leq p<\infty$. Let us prove that the group $\langle\omega_0 \rangle$ generated by $\omega_0 $ is dense in $\mathbb{T}^H_p$. Since the case $p=0$ was proved by S.~Rolewicz \cite{Rol}, we assume that $p>0$.

Let $\varepsilon>0$ and $\omega =(z_n)\in \mathbb{T}^H_p$, where $z_n ={\rm e}^{2\pi i y_n}, y_n \in [-\frac{1}{2}; \frac{1}{2})$. Set $q=\max(p,1)$ and choose $n$ such that
\begin{equation} \label{7}
\sum_{s=n}^{\infty} |1-z_s |^p +\frac{(n-1)2^p \pi^p}{2^{pn}} +\frac{1}{2^{pn}}\frac{4^p \pi^p}{2^p -1} <\varepsilon^q .
\end{equation}
By the definition of $\omega_0$, we can choose $k\leq k_n$ and integers $m_1,\dots, m_{n-1}$ such that
\[
|ka_s - y_s -m_s|<\frac{1}{2^n} \mbox{ for all } s=1,\dots, n-1.
\]
It is remained to prove that $d_p (\omega, k\omega_0 )< \varepsilon$.

For $s=1,\dots, n-1$, by (\ref{2}), we have
\begin{equation} \label{8}
| z_s - (z_s^0)^k |^p = |1 -\exp\{ 2\pi i (ka_s - y_s)\} |^p < 2^p \pi^p |ka_s - y_s -m_s|^p < \frac{2^p \pi^p}{2^{pn}}.
\end{equation}
For $s\geq n$  we have $(z_s^0)^k ={\rm e}^{2\pi i k a_s}$. Since $k\leq k_n$, by (\ref{2}), we obtain  \begin{equation} \label{9}
| z_s - (z_s^0)^k |^p \leq |1- z_s|^p + |1- (z_s^0)^k |^p < |1- z_s|^p +\frac{2^p \pi^p}{2^{ps}}.
\end{equation}
Then, by (\ref{7})-(\ref{9}), we have
\[
d_p^q (\omega, k\omega_0 )\leq \sum_{s=1}^{n-1} | z_s - (z_s^0)^k |^p +\sum_{s=n}^{\infty} \left( |1- z_s|^p +\frac{2^p \pi^p}{2^{ps}} \right) <
\]
\[
\frac{(n-1)2^p \pi^p}{2^{pn}} +\sum_{s=n}^{\infty} |1-z_n |^p +\frac{1}{2^{pn}}\frac{4^p \pi^p}{2^p -1} < \varepsilon^d.
\]
Thus $d_p (\omega, k\omega_0 )< \varepsilon$ and  $\langle\omega_0 \rangle$ is dense. $\Box$

By inequality (\ref{5}), we will consider the groups $\mathbb{T}^H_p, p=0$ or $1\leq p$, under the following metrics: if $\omega_j =(z_n^j)= ({\rm e}^{2\pi i \varphi_n^j})$, where $\varphi_n^j \in [-\frac{1}{2}; \frac{1}{2}), j=1,2$, then
\[
\rho_p (\omega_1,\omega_2)=\left( \sum_{n=1}^{\infty} |\varphi_n^1 -\varphi_n^2 |^p \right)^{1/p}, \mbox{ if } 1\leq p, \mbox{ and}
\]
\[
\rho_0 (\omega_1,\omega_2)= \sup \left\{ |\varphi_n^1 -\varphi_n^2 |, \; n=1,2,\dots \right\}, \mbox{ if } p=0.
\]
In the sequel we need some notations. For $p\geq 1$, nonnegative integers $k$ and $m$ and $\chi =(n_1, \dots, n_k, 0,\dots)\in \mathbb{Z}_0^{\infty}$, we set
\begin{itemize}
\item[-] $|\chi|_p := \sqrt[p]{|n_1|^p +\dots + |n_k|^p} $;
\item[-] $l(\chi)$ is the number of nonzero coordinates of $\chi$;
\item[-] $A(k,m) =\left\{ \chi =(0,\dots, 0, n_{m+1},\dots, n_s, 0,\dots)\in \mathbb{Z}_0^{\infty} : \;\; |\chi|_1 \leq k+1 \right\}$;
\item[-] the integral part of a real number $x$ is denoted by $[x]$.
\end{itemize}
We need the following lemma.
\begin{lem} \label{l1}
{\it For any $p\geq 1$ and $0<\varepsilon <1/4$ we set
\[
A_{\varepsilon,q} = \left\{ \chi\in\mathbb{Z}_0^{\infty} :\;\; 4\varepsilon |\chi|_q \leq 1 \right\}.
\]
Put $a=\left[ \frac{1}{4\varepsilon} \right] -1$ and $b= \left[ \left( \frac{1}{4\varepsilon} \right)^{q} \right]$. Then for every $p>1$ we have}
\[
A(a, 0) \subset A_{\varepsilon,p} \subset A(b,0).
\]
\end{lem}

\pr Let $\chi=(n_k) \in\mathbb{Z}_0^{\infty}$. Since $q\geq 1$ and $n_k \in \mathbb{Z}$, we have
\[
|\chi |_1 = \sum | n_k | \leq \sum | n_k |^q = |\chi |_q^q \leq \left( \sum | n_k | \right)^q =|\chi |_1^q
\]
and the assertion follows. $\Box$

Let us consider the sequence $\mathbf{e} =\{ e_n\} \in\mathbb{Z}_0^{\infty}$, where $e_1 =(1,0,0,\dots), e_2 = (0,1,0,\dots), \dots$. Denote by $(\mathbb{Z}_0^{\infty} , \mathbf{e})$ the group $\mathbb{Z}_0^{\infty}$ equipped with the finest Hausdorff group topology for which $e_n$ converges to zero (in fact, it is a Graev free topological abelian group over the convergent sequence $\mathbf{e}\cup \{ 0\}$).

\begin{teo} \label{t1}
{\it
\begin{enumerate}
\item $\mathbb{T}^H_p$ is not locally quasi-convex and, so, not reflexive for any $1< p <\infty$.
\item $\mathbb{T}^H_1$ is reflexive and, hence, locally quasi-convex.
\item $\mathbb{T}^H_0$ is reflexive and, hence, locally quasi-convex. It is not a {\rm QI}-group.
\item $\left(\mathbb{T}^H_p \right)^{\wedge}$ is topologically isomorphic to $\left(\mathbb{T}^H_0 \right)^{\wedge}$ and, hence, reflexive for any $1< p <\infty$.
\item $\left(\mathbb{T}^H_0 \right)^{\wedge}= (\mathbb{Z}_0^{\infty} , \mathbf{e})$.
\item $\left(\mathbb{T}^H_1 \right)^{\wedge}$ is algebraically isomorphic to $\mathbb{Z}^{\infty}_b$.
\end{enumerate} }
\end{teo}

\pr {\bf 1.} Evidently, $\mathbb{T}^n , n\geq 1,$ are  closed
subgroups of $\mathbb{T}^H_p$. Let $\chi\in \left(\mathbb{T}^H_p\right)^{\wedge}, 1\leq p<\infty$.
Then $\chi |_{\mathbb{T}^n}$ is a character of $\mathbb{T}^n$.
Hence $\chi |_{\mathbb{T}^n} = (m_1 , \dots , m_n ), m_n \in
\mathbb{Z}$.

a) {\it Let us prove that $\left(\mathbb{T}^H_p \right)^{\wedge}$ is algebraically isomorphic to $\mathbb{Z}^{\infty}_0$ for every $1<p<\infty$}.

It is clear that $\mathbb{Z}^{\infty}_0 \subset \left(\mathbb{T}^H_p \right)^{\wedge}$.
For the converse inclusion it is remained to prove that only finite number of integers $m_n$
are nonzero. Assume the converse and $m_{s_l } \not= 0, l=1,2,\dots$. We can assume that $m_{s_l} >0$. Set
\[
a_1 =\frac{1}{2\pi} , \; a_2 =\frac{-1}{2\pi\cdot 2\ln 2} , \; a_3
=\frac{-1}{2\pi\cdot 3\ln 3} , \dots , a_{k_1}
=\frac{-1}{2\pi\cdot k_1\ln k_1} ,
\]
where $k_1$ is the first number  such that $\sum_{k=1}^{k_1} a_k <
\frac{-1}{2\pi}$,
\[
a_{k_1 +1} =\frac{1}{2\pi\cdot (k_1 +1)\ln (k_1 +1)} , \dots , a_{k_2}
=\frac{1}{2\pi\cdot k_2\ln k_2} ,
\]
where $k_2$ is the first number  such that $\sum_{k=1}^{k_2} a_k >
\frac{1}{2\pi}$, and etc. Put $z_{s_l} =\exp (2\pi i\frac{a_l}{m_{s_l}})$ and $z_n =1$
for the remainder $n$. Obviously, $\omega =( z_n ) \in \mathbb{T}^H_p$.
Since $\cup_n \mathbb{T}^n$ is dense in $\mathbb{T}^H_p$ and $\chi$
is continuous, there exists $(\chi , \omega ) = \lim_l \exp
(2\pi i\sum_{n=1}^{k_l} a_n )$. But $\textrm{Im} \exp
(i\sum_{n=1}^{k_l} a_n )
> \sin 1,$ for even $l$, and $< - \sin 1$, for odd $l$. It is a
contradiction.

b) {\it Let us prove that $\left(\mathbb{T}^H_1 \right)^{\wedge}$ is algebraically isomorphic to $\mathbb{Z}^{\infty}_b$}.

For the inclusion $\left(\mathbb{T}^H_1 \right)^{\wedge} \subset \mathbb{Z}^{\infty}_b$ we need to prove that $\{ m_k \}$ is bounded. Assuming the converse, we can choose a subsequence $k_l$ such that  $|m_{k_l}|> l^2, l=1,2,\dots $.
Set $\omega =(z_n)$, where $z_{n} = \exp (2\pi i /3m_{k_l} )$, if $n=k_l$, and $z_n =1$ otherwise. It is clear that $\omega\in \mathbb{T}^H_1$ and $(\chi ,\omega)$ does not exist.

On the other hand, if $|\textbf{n}|_{\infty} <\infty$, then $\chi=\textbf{n}$ is a continuous character of $\mathbb{T}^H_1$. Thus $\mathbb{Z}^{\infty}_b \subset \left(\mathbb{T}^H_1 \right)^{\wedge}$.

{\bf 2.} Set $U_{\varepsilon}$ is the $\varepsilon$-neighborhood of the unit in $\mathbb{T}^H_p , 1\leq p<\infty$.

a) {\it Let $1<p<\infty$. We will prove that $A_{\varepsilon ,q} = U_{\varepsilon}^{\triangleright}$ for any $0< \varepsilon <1$.} By Lemma \ref{l1}, this shows that the sets $A(k,0)$ form a decreasing family of precompact sets such that any compact set in the hemicompact group $\left(\mathbb{T}^H_p\right)^{\wedge}$ is contained in some $A(k,0)$.

Let $\chi \in U_{\varepsilon}^{\triangleright}$. If $\chi =(n_1 , \dots, n_k , 0,\dots )$ and $\omega =(z_j)=({\rm e}^{2i\pi \varphi_j})$, where $\varphi_j  \in [-\frac{1}{2}; \frac{1}{2})$, then
\[
(\chi , \omega)= z_1^{n_1} \dots z_k^{n_k} = \exp ( 2i \pi (n_1 \varphi_1 + \dots +n_k \varphi_k )).
\]
Since $U_{\varepsilon}$ is pathwise connected, then
\begin{equation} \label{12}
{\rm Re}(\chi , \omega) \geq 0 , \forall \omega\in U_{\varepsilon}, \mbox{ iff } -\frac{1}{4} \leq n_1 \varphi_1 + \dots +n_k \varphi_k \leq \frac{1}{4}, \forall \omega\in U_{\varepsilon},
\end{equation}
and, in particular, for all
\begin{equation} \label{13}
\omega_k = (z_1 ,\dots, z_k)\in \mathbb{T}^k \mbox{ such that } \sum_{j=1}^k | \varphi_j |^p < \varepsilon^p .
\end{equation}
By H\"{o}lder's inequality, we have
\[
| n_{1} \varphi_{1} +\dots +n_{k} \varphi_{k} | \leq |\chi |_q \cdot \sqrt[p]{ |\varphi_{1}|^p +\dots +|\varphi_{k}|^p } \leq \varepsilon |\chi|_q.
\]

It is clear that the supremum  of the function $f=n_1 \varphi_1 + \dots +n_k \varphi_k$ under the condition (\ref{13}) is achieved when $\varphi_i =\varepsilon\cdot {\rm sign} (n_i) \left(\frac{|n_i|}{|\chi |_q} \right)^{q/p}$ and equals to $\varepsilon |\chi|_q$.  By (\ref{12}), we have
\[
U_{\varepsilon}^{\triangleright} = \left\{ \chi \in \left(\mathbb{T}^H_p\right)^{\wedge} : \; \sup f (\omega_k) \leq \frac{1}{4}, \omega \in U_{\varepsilon} \right\} = \left\{ \chi \in \mathbb{Z}_0^{\infty} :\; 4 \varepsilon |\chi|_q \leq 1 \right\} .
\]

b) {\it Let  $p=1$ and $0<\varepsilon <\frac{1}{2}$. Set $Z_{\varepsilon} = \left\{ \mathbf{n} \in \mathbb{Z}^{\infty}_b : \; | \mathbf{n} |_b \leq \frac{1}{\varepsilon} \right\}$. We will prove that}
\[
Z_{4\varepsilon} \subset U_{\varepsilon}^{\triangleright} \subset Z_{2\varepsilon}.
\]
This shows that the sets $Z_{\varepsilon}$ form a decreasing family of precompact sets and each compact set in the hemicompact group $\left(\mathbb{T}^H_1\right)^{\wedge}$ is contained in some $Z_{\varepsilon}$.

If $\omega\in U_{\varepsilon}$, then analogously to case a), we have the following: since $U_{\varepsilon}$ is pathwise connected, then
\begin{equation} \label{14}
\chi \in U_{\varepsilon}^{\triangleright} \Leftrightarrow |\sum n_i \varphi_i | \leq \frac{1}{4} ,\; \forall \omega \in U_{\varepsilon}.
\end{equation}
Since
\[
|\sum n_i \varphi_i | \leq |\mathbf{n}|_b \cdot \sum |\varphi_i| = |\mathbf{n}|_b \cdot |\omega |_1 \leq |\mathbf{n}|_b \cdot \varepsilon,
\]
we have
\[
Z_{4\varepsilon} = \left\{ \mathbf{n} \in \left(\mathbb{T}^H_1\right)^{\wedge} : |\mathbf{n}|_b \leq \frac{1}{4\varepsilon} \right\} \subset U_{\varepsilon}^{\triangleright}.
\]
Let us prove the second inclusion. If $|\mathbf{n}|_b \geq \frac{1}{2\varepsilon}$, then $|n_j|\geq \frac{1}{2\varepsilon}$ for some $j$. If $\omega=(z_n)$, where $z_n = {\rm e}^{2\pi i\varepsilon \frac{3}{4}}$ if $n=j$ and $z_n =1$ otherwise, then $\omega\in U_{\varepsilon}$ and  $|\sum n_i \varphi_i | =|n_j| \frac{3}{4}\varepsilon \geq 3/8$. This contradicts to (\ref{14}). Thus $U_{\varepsilon}^{\triangleright} \subset Z_{2\varepsilon}$.

{\bf 3.} a) {\it Let $1<p<\infty$. Let us show that for each neighborhood $W$ of $\chi =0$ in $(\mathbb{T}_p^H)^{\wedge}$ and a positive integer $k$ there exists $m$ such that}
\begin{equation} \label{15}
A(k,m) \subset W.
\end{equation}

According to corollary  4.4 \cite{Aus}, there exists a sequence $\{ a_n \}$ such that $a_n \to 1$ and $\{ a_n \}^{\vartriangleright} \subset W$.  Let us show that
$A(k,m) \subset \{ a_n \}^{\vartriangleright}$ for all large $m$. Set $a_n = (z^n_k) = ({\rm e}^{2\pi i\varphi^n_k} ) $.

Let $q$ be such that $\frac{1}{p} +\frac{1}{q}=1$. Set $A=\max \{ \sqrt[q]{|n_1|^q +\dots +|n_{k+1}|^q} ,\; |n_j|\leq k+1\}$. Thus, if $\chi\in A(k,0)$, then $l(\chi)\leq k+1$ and $|\chi|_q \leq A$.
Now choose $N_1$ such that $\rho_p (e, a_n ) \leq \frac{1}{4A}, \forall n\geq N_1 ,$ and choose $M>N_1$ such that
\[
\sum_{k=M+1}^{\infty} |\varphi^n_k|^p \leq \frac{1}{(4A)^p} , \forall n \leq N_1.
\]
Then for all $n$ we have
$\sum_{k=M+1}^{\infty} |\varphi^n_k|^p < \frac{1}{(4A)^p} .$
Hence, for all $\chi\in A(k,m), m\geq M,$ and $a_n$, by H\"{o}lder's inequality, we have
\[
| n_{m+1} \varphi_{m+1}^n +\dots +n_{l} \varphi_{l}^n | \leq |\chi |_q \cdot \sqrt[p]{ |\varphi_{m+1}^n|^p +\dots +|\varphi_{l}^n|^p } \leq A \cdot \frac{1}{4A} =\frac{1}{4} .
\]
Therefore $A(k,m) \subset \{ a_n \}^{\vartriangleright}$ for all $m\geq M$.

b) {\it Let $p=1, \varepsilon >0$ and $W$ be an open neighborhood of the neutral element of $\left(\mathbb{T}^H_1\right)^{\wedge}$. Set
\[
Z_{\varepsilon}^l := \left\{ \mathbf{n} =(0,\dots , 0, n_{l+1}, n_{l+2},\dots ) : \; |\mathbf{n} |_b \leq \frac{1}{\varepsilon} \right\} \subset Z_{\varepsilon}.
\]
Let us show that $Z_{\varepsilon}^l \subset  W$ for all large $l$.}

Analogously, according to corollary  4.4 \cite{Aus},  there exists a sequence $\{ a_k \}$ such that $a_k \to 1$ and $\{ a_k \}^{\vartriangleright} \subset W$.  It is enough to show that
$Z_{\varepsilon}^l \subset \{ a_k \}^{\vartriangleright}$ for all large $l$.

Choose $N$ such that $\rho_1 (e, a_k ) < \frac{\varepsilon}{4}, \forall k> N$. Choose $l_0 >N$ such that
\[
\sum_{i=l_0+1}^{\infty} |\varphi^k_i | < \frac{\varepsilon}{4} , \; \mbox{ for every } k=1,\dots,N.
\]
Then the last inequality is true for all $k$. Therefore for every $l\geq l_0,$ every  $\chi \in Z_{\varepsilon}^l$ and $a_k$, we have
\[
| \sum_{i=1}^{\infty} n_i \varphi_i^k | =| \sum_{i=l+1}^{\infty} n_i \varphi_i^k | \leq |\mathbf{n} |_b \sum_{i=l+1}^{\infty} |\varphi_i^k| < \frac{1}{\varepsilon} \cdot \frac{\varepsilon}{4} = \frac{1}{4} .
\]
Hence  $\mathbf{n} \in \{ a_n\}^{\triangleright}$. Thus $Z_{\varepsilon}^l \subset \{ a_n \}^{\vartriangleright}$.

{\bf 4.} a) {\it Let $1<p<\infty$ and $\varepsilon <0,01$. Let $\chi_{\alpha} \to \chi$, where $\chi=(n_1,\dots, n_s, 0,\dots), \chi_{\alpha} , \chi \in U_{\varepsilon}^{\triangleright}$. Let us prove that for every $M$ there exists $\alpha_0$ such that}
\[
\chi_{\alpha} = (n_1,\dots, n_s, 0,\dots, 0_M, n_{M+1}^{\alpha},\dots) , \forall \alpha \geq \alpha_0.
\]

By item 2a and Lemma \ref{l1}, we have $U_{\varepsilon}^{\triangleright} = A_{\varepsilon ,q} \subset A\left( \left[ \left( \frac{1}{4\varepsilon} \right)^{q} \right] , 0\right)$. Thus
\begin{equation} \label{16}
| n_j |+ |n_j^{\alpha} | \leq 2\left[ \left( \frac{1}{4\varepsilon} \right)^{q} \right] +2.
\end{equation}
Set $q=\max\left\{ \left( 2\left[ \left( \frac{1}{4\varepsilon} \right)^{q} \right] +2\right)^2 +1 , \frac{2}{\varepsilon^2} \right\}$. Let a sequence $\{ a_n \}$ be such that $a_n \to e$ and consists the following elements
\[
(\exp \left( 2\pi i\frac{k_1}{q} \right), \dots , \exp \left( 2\pi i\frac{k_M}{q} \right), 1, \dots ), \mbox{ where } \; k_i =0, \pm 1,\dots , \pm q.
\]
Since $\{ a_n \}$ is compact, $\{ a_n \}^{\triangleright}$ is open. Thus there exists $\alpha_0$ such that  $\chi -\chi_{\alpha} \in \{ a_n \}^{\triangleright}$ for $\alpha >\alpha_0$. In particular,
\begin{equation} \label{17}
{\rm Re} \left\{ \exp \left(2i\pi (n_j - n_j^{\alpha} ) \frac{k}{q} \right)\right\} \geq 0, \; j=1, \dots , M, |k| \leq q.
\end{equation}
Now we assume the converse and $| n_j - n_j^{\alpha} |>0$ for some $0<j\leq M$.  Then, by (\ref{16}), $1\leq | n_j - n_j^{\alpha} | \leq | n_j |+ |n_j^{\alpha} | \leq \sqrt{q}$. Since $q>\frac{2}{\varepsilon^2}$ and $\varepsilon < 0,01$, then
\[
\frac{1}{q} \leq \frac{| n_j - n_j^{\alpha} | }{q} \leq \frac{\sqrt{q}}{q} < \varepsilon <0,01.
\]
Hence there exists $|k_j |\geq 1$ such that $\frac{1}{4} < \frac{( n_j - n_j^{\alpha} ) k_j}{q} <\frac{1}{2}$. Therefore for this $k_j$ the inequality (\ref{17}) is wrong.

b) {\it Let us prove that $A(k,m)$ is compact $(1<p<\infty)$.}

Let a net $\{\chi_{\alpha}\} \subset A(k,m)$ is fundamental. Since, by item 2, $A(k,m)$ is precompact and is contained in some $U_{\varepsilon}^{\triangleright}$, then it converges to some $\chi =(n_1,\dots, n_s, 0,\dots)$. Choose $\alpha_0$ such that $\chi_{\alpha} =(n_1,\dots, n_s, 0,\dots, 0_M, n_{M+1}^{\alpha},\dots) $, $\forall \alpha \geq \alpha_0$. Since $|\chi |_1 \leq |\chi_{\alpha} |_1 \leq k+1$, we obtain that $\chi \in A(k,m)$.

{\bf 5.} {\it Let us prove that $\left( \mathbb{T}^H_p \right)^{\wedge\wedge} = \mathbb{T}^{\infty}_0, 1<p<\infty$.}

It is clear that $\left( \mathbb{T}^H_p \right)^{\wedge\wedge} \subset \left( (\mathbb{Z}^{\infty}_0 )_d \right)^{\wedge} =\mathbb{T}^{\infty}$.
Let $\omega =(z_n)=({\rm e}^{2i\pi \varphi_n} ) \in \left( \mathbb{T}^H_p \right)^{\wedge\wedge}$.

a) {\it Let us show that $\left( \mathbb{T}^H_p \right)^{\wedge\wedge} \subset \mathbb{T}^{\infty}_0$, i.e. $z_n\to 1$.}

Assume the converse and $z_n\not\to 1$. Then there exists a subsequence $n_k$ such that $\varphi_{n_k} \to \alpha\not= 0$. We will show that $\omega$ is discontinuous at $0$. Let $W$ be a neighborhood of the neutral element of $(\mathbb{T}_p^H)^{\wedge}$. By (\ref{15}), there exists $m$ such that $A(1,m) \subset W$. In particular, if $n_k >m$, then $\chi_k =(0,\dots, 0,1,0,\dots)$, where $1$ occupies position $n_k$, belongs to $W$. Then
\[
(\omega , \chi_k ) = \exp (2i\pi \varphi_{n_k} ) \to \exp(2i\pi \alpha) \not= 1.
\]
Thus $\omega$ is discontinuous. Hence $\left( \mathbb{T}^H_p \right)^{\wedge\wedge} \subset \mathbb{T}^{\infty}_0$.

b) {\it Let us prove the converse inclusion: $\left( \mathbb{T}^H_p \right)^{\wedge\wedge} \supset \mathbb{T}^{\infty}_0$, i.e. if $z_n \to 1$, then $\omega$ is a continuous character of $\left( \mathbb{T}^H_p \right)^{\wedge}$.}

Since $(\mathbb{T}^H_p)^{\wedge}$ is a hemicompact $k$-space \cite{Cha}, by item 2a, it is enough to prove that $\omega$ is continuous on $A(k,0)$. Let $\varepsilon >0$ and $\chi_{\alpha} \to \chi$, where $\chi_{\alpha} , \chi \in A(k,0)$. Let $q$ be such that $\frac{1}{p} +\frac{1}{q}=1$. Set $A=\max \{ \sqrt[q]{|n_1|^q +\dots +|n_{k+1}|^q} ,\; |n_j|\leq k+1\}$. Thus, if $\eta\in A(k,0)$, then $|\eta|_q \leq A$.

Choose $M$ such that $\sqrt[p]{\sum_{k=1} \left|\varphi_{M+k} \right|^p } < \frac{\varepsilon}{2\pi A}$. By item 4a, for $M$ we can choose $\alpha_0$ such that
\[
\chi_{\alpha} = (n_1,\dots, n_s, 0,\dots, 0_M, n_{M+1}^{\alpha},\dots), \mbox{ where } \chi=(n_1,\dots, n_s, 0,\dots), \forall \alpha >\alpha_0.
\]
Then $
(\omega , \chi_{\alpha} -\chi) = \exp \left( 2i\pi \sum_{k=1}^{\infty} n_{M+k}^{\alpha} \varphi_{M+k} \right) .$ By H\"{o}lder's inequality, we have
\[
| \sum n_{M+k}^{\alpha} \varphi_{M+k} |\leq |\chi_{\alpha} |_q  \cdot \sqrt[p]{ \sum \left|\varphi_{M+k} \right|^p }< A\cdot \frac{\varepsilon}{2\pi A} =\frac{\varepsilon}{2\pi}.
\]
Therefore for $\alpha >\alpha_0$, by (\ref{2}), we obtain
\[
| (\omega , \chi_{\alpha}) -(\omega , \chi) | = | 1- (\omega , \chi_{\alpha} -\chi) | <\varepsilon
\]
and $\omega$ is continuous.

c) {\it Let us prove that $(\mathbb{T}^H_p)^{\wedge\wedge}$ is topologically isomorphic to $\mathbb{T}^H_0$.}

Since $\mathbb{T}^H_p$ is Polish, then $(\mathbb{T}^H_p)^{\wedge\wedge}$ is Polish too \cite{Cha}. Since $(\mathbb{T}^H_p)^{\wedge\wedge}$ and $\mathbb{T}^H_0$ are the same Borel subgroup of $\mathbb{T}^{\infty}$, they must coincide topologically.

{\bf 6.} a) {\it Let us prove that $\mathbb{Z}^{\infty}_0$ is dense in $\left( \mathbb{T}_1^H\right)^{\wedge} =\mathbb{Z}^{\infty}_b$.}

Let $\mathbf{n}_0 =(n_i) \in \mathbb{Z}^{\infty}_b$ and $W$ be a neighborhood of the neutral element. Let $\varepsilon$ be such that $|\mathbf{n}_0 |_b \leq \frac{1}{\varepsilon}$. By item 3b of the proof, we can choose $l$ such that $Z_{\varepsilon }^l \subset W$. Set $\mathbf{n} = (n_1, \dots,n_l,0,\dots) \in \mathbb{Z}^{\infty}_0$. Then $\mathbf{n}_0 -\mathbf{n} \in Z^{l}_{\varepsilon} \subset W$ q.e.d.

b) {\it Let us prove that $\mathbb{T}^H_1$ is reflexive.}

Set $t : \mathbb{T}^H_1 \to \mathbb{T}^{\infty}$ is the natural continuous monomorphism. Since the image of $t$ is dense, $t^{*} : (\mathbb{T}^{\infty})^{\wedge} =\mathbb{Z}^{\infty}_0 \to (\mathbb{T}^H_1)^{\wedge} =\mathbb{Z}^{\infty}_b$ is injective. As it was proved in part a), $t^{*}$ has the dense image. Hence $t^{**} : (\mathbb{T}^H_1)^{\wedge\wedge} \to \mathbb{T}^{\infty}$ is injective. By corollary 3 \cite{Cha},
it is enough to prove that $\mathbb{T}^H_1 =(\mathbb{T}^H_1)^{\wedge\wedge}$ algebraically.
Let $\omega =(z_n)=({\rm e}^{2\pi i \varphi_n}) \in (\mathbb{T}^H_1)^{\wedge\wedge}$. If $\sum |\varphi_n | =\infty$, then, in the standard way, we can construct $\mathbf{n} =(\pm 1)$ such that $(\omega , \mathbf{n})$ does not exist. Thus $\omega$ must be contained in $\mathbb{T}^H_1$.

{\bf 7.} a) {\it  Let $1<p<\infty$. Let us prove that
\begin{enumerate}
\item $\mathbb{T}^H_p$ is not locally quasi convex.
\item $(\mathbb{T}^H_p)^{\wedge}$ is reflexive.
\item $(\mathbb{T}^H_p)^{\wedge}= (\mathbb{T}^H_0)^{\wedge}$ and, hence,  $(\mathbb{T}^H_p)^{\wedge}$ does not depend on $p$.
\item $\mathbb{T}^H_0$ is reflexive.
\end{enumerate} }

Let $\alpha_p : \mathbb{T}^H_p \to (\mathbb{T}^H_p)^{\wedge\wedge} = \mathbb{T}^H_0$ be the canonical homomorphism. Then $\alpha_p$ has the dense image. By Proposition \ref{p3}.1, $(\mathbb{T}^H_p)^{\wedge}$ and $\mathbb{T}^H_0 = (\mathbb{T}^H_p)^{\wedge\wedge}$ are reflexive. Thus, $(\mathbb{T}^H_p)^{\wedge}=(\mathbb{T}^H_0)^{\wedge}$ does not depend on $p$. By Proposition \ref{p3}.2, $\mathbb{T}^H_p$ is not locally quasi convex.

b) {\it Let us prove that $\mathbb{T}^H_0$ is not a QI-group.}

Assume the converse and $\mathbb{T}^H_0$ is a QI-group.
Assume that $\mathbb{T}^H_0$ is represented in $\mathbb{T}^{\infty}$. Set
\[
U_{\varepsilon}^d =\{ h\in \mathbb{T}^H_0 :\; d(e,h)<\varepsilon \} , \;\; U_{\varepsilon} =\{ h\in \mathbb{T}^H_0 :\; \rho_0 (e,h)<\varepsilon \}.
\]
By Proposition \ref{p1}, there exists $\varepsilon_0 >0$ such that ${\rm Cl}_{\mathbb{T}^{\infty}} (U_{\varepsilon_0}^d)$ is compact in $\mathbb{T}^{\infty}$ and it is contained in $\mathbb{T}^H_0$. Since the Polish group topology is unique, there exists $\varepsilon >0$, such that $U_{\varepsilon} \subset U_{\varepsilon_0}^d$.
Set $z=\exp \left( 2\pi i \frac{\varepsilon}{2\pi}\right)$. Then $\omega_k = (z,\dots , z, 1_{k+1}, 1,\dots) \in {\rm Cl}_{\mathbb{T}^{\infty}} (U_{\varepsilon}) \subset \mathbb{T}^H_0$ for every $k$.  But $\omega_k$ converges in $\mathbb{T}^{\infty}$ to $\omega = (z) \not\in \mathbb{T}^H_0$. Hence $\mathbb{T}^H_0$ can not be represented in $\mathbb{T}^{\infty}$. If $\mathbb{T}^H_0$ is represented in another locally compact group $X$, then, by 25.31(b) \cite{HR}, $X = \mathbb{T}^{\infty}\times X_1$ and $\mathbb{T}^H_0$ condensates to $\mathbb{T}^{\infty}$. Hence $\mathbb{T}^H_0$ is not a QI-group.

{\bf 8.} {\it Let us prove that $\left(\mathbb{T}^H_p \right)^{\wedge}$ and $(\mathbb{Z}_0^{\infty} , \mathbf{e})$ are topologically isomorphic.}

Since $e_n \in A(0,n)$ for every $n$, $e_n$ converges to zero in $\left(\mathbb{T}^H_p \right)^{\wedge}$ by (\ref{15}). Thus, by definition, $id: (\mathbb{Z}_0^{\infty} , \mathbf{e}) \mapsto \left(\mathbb{T}^H_p \right)^{\wedge}$ is continuous.

Let us prove the  converse, i.e. if $W$ is an open neighborhood of zero in $(\mathbb{Z}_0^{\infty} , \mathbf{e})$, then $W$ is open in $\left(\mathbb{T}^H_p \right)^{\wedge}$. At first we note that since $\left(\mathbb{T}^H_p \right)^{\wedge}$ is a hemicompact $k$-space, a set $Y$ is closed if and only if $Y\cap A(k,0)$ is closed in $A(k,0)$. By the construction of $(\mathbb{Z}_0^{\infty} , \mathbf{e})$ \cite{ZP1}, we may assume that $W$ has the form
\[
W=\cup_{k=1}^{\infty} (A^{\ast}_{i_1} +A^{\ast}_{i_2} +\dots + A^{\ast}_{i_k}),
\]
 where  $1\leq i_1 <i_2<\dots, A^{\ast}_{n} =\{ 0, \pm e_m :\; m\geq n \} .$
Now assume the converse and $W$ is not open in $\left(\mathbb{T}^H_p \right)^{\wedge}$. Then there exists $k_0$ such that $A(k_0, 0)\setminus W$ is not closed in $A(k_0, 0)$. In particular, $A(k_0, 0)\setminus W \not= \emptyset$. Let a net $\{ \chi_{\alpha}\} \subset A(k_0, 0)\setminus W$ be such that $\chi_{\alpha}$ converges to $\chi_0 =(n_1,\dots, n_s, 0,\dots)\in A(k_0, 0)\cap W$. Since $(\mathbb{Z}_0^{\infty} , \mathbf{e})$ is a topological group, we can find a neighborhood $W_1$ of zero such that
\[
\chi_0 + W_1 \subset W, \mbox{ where } W_1 =\cup_{k=1}^{\infty} (A^{\ast}_{j_1} +A^{\ast}_{j_2} +\dots + A^{\ast}_{j_k}), 1\leq j_1 <j_2<\dots
\]
By the construction of $W_1$, we have $A(k_0 , j_{k_0 +1}) \subset W_1$. By 4a, we can choose $\alpha_0$ such that
\[
\chi_{\alpha} = (n_1,\dots, n_s, 0,\dots, 0_{j_{k_0 +1}}, n_{j_{k_0 +1}+1}^{\alpha},\dots)\in A(k_0 ,0) , \forall \alpha \geq \alpha_0.
\]
Thus $\chi_{\alpha} -\chi_0 \in A(k_0 , j_{k_0 +1})$ for every $\alpha \geq \alpha_0$. So
\[
\chi_{\alpha} = \chi_0 + (\chi_{\alpha} -\chi_0) \in \chi_0 +W_1 \subset W, \; \forall \alpha \geq \alpha_0.
\]
It is impossible since $\chi_{\alpha}\not\in W$. Thus $W$ is open in $\left(\mathbb{T}^H_p \right)^{\wedge}$. The theorem is proved. $\Box$

\begin{rem}
Let $G$ be the group $\mathbb{Z}_0^\infty$ with the discrete topology and $H=(\mathbb{Z}_0^\infty , \mathbf{e})$. Then $G$ and $H$ are reflexive groups and $id : G\to H$ is a continuous isomorphism, but
$id^{\ast} : H^{\wedge}=\mathbb{T}_0^H \to G^{\wedge}= \mathbb{T}^\infty$ is only a continuous injection. $\Box$
\end{rem}

\begin{rem}
We can prove that $(\mathbb{T}_0^H)^\wedge =\mathbb{Z}_0^\infty$ algebraically using the following observation.

{\it Let $X, G,$ and $Y$ be a topological abelian groups such that there exist continuous monomorphisms $i: X\to G$ and $j: G\to Y$ with the dense image. If $i^{\ast}\circ j^{\ast}$ is bijective, then $i^{\ast}$ and $j^{\ast}$ are bijective too.}

(Since $i$ has the dense image, $i^{\ast}$ is injective. Since $i^{\ast}\circ j^{\ast}$ is bijective, then $i^{\ast}$ must be surjective. Hence $i^{\ast}$ is bijective. Thus $j^{\ast} = (i^{\ast})^{-1} (i^{\ast}\circ j^{\ast})$ is bijective.)

Now, since $\mathbb{T}^H_p \subset \mathbb{T}^H_0 \subset \mathbb{T}^{\infty}$ and $(\mathbb{T}^H_p )^{\wedge} =(\mathbb{T}^{\infty})^\wedge = \mathbb{Z}^{\infty}_0$, then, by the observation, $(\mathbb{T}^H_0 )^{\wedge} =\mathbb{Z}^{\infty}_0$ algebraically. $\Box$
\end{rem}

Now we consider the group $G_p$, $1<p<\infty$.
Since $Q$ is dense in $G_p$, $G_p^{\wedge} \subset Q_d^{\wedge} = \Delta_{\mathbf{a}}$, where $\mathbf{a} =(a_1, a_2,\dots)$ (here $Q_d$ denotes the group $Q$ with discrete topology). By section 25.2 \cite{HR}, we have
\[
(\chi , z) = z^{\sum_{k=1}^{\infty} \omega_k \gamma (k)} , \, \forall z\in Q, \mbox{ where } \chi =(\omega_k)\in \Delta_{\mathbf{a}} , \omega_k \in \{ 0,1,\dots, a_k -1\} .
\]
We need the following lemma.
\begin{lem} \label{l2}
{\it The group $Q$ is dense in $G_0$.}
\end{lem}

\pr We identify $\mathbb{T}$ with $[0,1)$ and denote by $\langle x\rangle$ the distance of $x$ from the nearest integer. By (\ref{2}), we can consider the following equivalent metric on $G_0$
\[
r_0 (x_1,x_2) = \sup \{ \langle\gamma (n) (x_1 -x_2 )\rangle,\; n\in \mathbb{N} \}.
\]
If $x\in [0,1)$, we can write
\[
x=\sum_{k=1}^{\infty} \frac{\varepsilon_k (x)}{a_1 a_2\dots a_k} , \mbox{ where } \varepsilon_k (x) =0,\dots, a_k -1,
\]
and for every $M>0$ there exists $k>M$ such that $\varepsilon_k (x) <a_k -1$. Thus
\[
\gamma (n) x ({\rm mod 1}) = \sum_{k=n}^{\infty} \frac{\varepsilon_k (x)}{a_n \dots a_k} = \frac{\varepsilon_n (x)}{a_n} + \frac{\theta_n}{a_n},
\]
where  $0\leq \theta_n <1$, since
\[
 \theta_n = \frac{\varepsilon_{n+1}}{a_{n+1}} + \frac{\varepsilon_{n+2}}{a_{n+1} a_{n+2}}+
 \frac{\varepsilon_{n+3}}{a_{n+1} a_{n+2} a_{n+3}} +\dots <
 \frac{a_{n+1}-1}{a_{n+1}} + \frac{a_{n+2}-1}{a_{n+1}a_{n+2}}+ \frac{a_{n+3}-1}{a_{n+1}a_{n+2} a_{n+3}} +\dots
 \]
 \[
 = \left( 1 -\frac{1}{a_{n+1}} \right) + \left( \frac{1}{a_{n+1}} - \frac{1}{a_{n+1}a_{n+2}}\right) +
 \left(  \frac{1}{a_{n+1}a_{n+2}} - \frac{1}{a_{n+1}a_{n+2} a_{n+3}}\right) + \dots =1.
 \]

By the definition of $G_0$, we have
\begin{equation} \label{18}
\langle \gamma (n) x \rangle =\langle \frac{\varepsilon_n (x) +\theta_n}{a_n} \rangle \to 0.
\end{equation}
Now we set $x_N = \sum_{k=1}^{N-1} \frac{\varepsilon_k (x)}{a_1 a_2\dots a_k} \in Q$. Then
\[
x-x_N = \sum_{k=N}^{\infty} \frac{\varepsilon_k (x)}{a_1 a_2\dots a_k} = \frac{1}{a_1 a_2\dots a_{N-1}} \cdot \frac{\varepsilon_N (x) +\theta_N}{a_N}
\]
and
\begin{equation} \label{19}
\gamma (n) (x-x_N) = \left\{
\begin{array}{ll}
\gamma (n) (x)({\rm mod 1}), & \mbox{ for } N\leq n \\
\frac{1}{a_n \dots a_{N-1}} \cdot \frac{\varepsilon_N (x) +\theta_N}{a_N} , & \mbox{ for } 1\leq n <N
\end{array}
\right.
\end{equation}

Let $\varepsilon>0$. Since $a_n\to\infty$, by (\ref{18}), we can choose $N$ such that $\langle \gamma (n) x \rangle <\varepsilon$ for all $n\geq N$ and $a_{N-1} > 1/\varepsilon$. Then, by  (\ref{19}), we obtain $r_0 (x, x_N)<\varepsilon$. Thus $Q$ is dense in $G_0$. $\Box$

 For $1<p<\infty$ and $p=0$  we consider the following homomorphisms
\[
S_p : \, G_p\mapsto \mathbb{T}\times \mathbb{T}_p^H \thickapprox \mathbb{T}_p^H ,\quad S_p (z) =(z, z^{\gamma (2)}, \dots, z^{\gamma (k)},\dots).
\]
It is clear that $S_p$ is a topological isomorphism of $G_p$ onto the following closed subgroup of $\mathbb{T}_p^H$
\[
\{ \omega \in \mathbb{T}_p^H : \, \omega = (z^{\gamma (1)}, z^{\gamma (2)},\dots , z^{\gamma (n)},\dots ) \} .
\]
We will identify $G_p$ with this subgroup. Set $\boldsymbol{\gamma} =\{ \gamma_n \}$. Then, by Lemma \ref{l2}, $\boldsymbol{\gamma}$ is a $TB$-sequence.  Denote by $(\mathbb{Z},\boldsymbol{\gamma})$ the group of integers equipped with the finest Hausdorff group topology in which $\gamma_n $ converges to zero.

\begin{teo} \label{t2}
{\it  Let $1<p<\infty$. Then
\begin{enumerate}
\item $G_p$ is not locally quasi-convex and, hence, not reflexive.
\item $G_0$ is reflexive and, hence, locally quasi-convex.
\item $\left( G_p \right)^{\wedge}=\left( G_0 \right)^{\wedge} =(\mathbb{Z},\boldsymbol{\gamma})$ and, hence, reflexive.
\end{enumerate} }
\end{teo}

\pr {\bf 1.} {\it Let us prove that $G_p , 1<p<\infty$ or $p=0,$ is dually closed in
$\mathbb{T}_p^H$.}

Let $\omega =(z^{\gamma (n)})$ and $\omega_0 = (z_1 , z_2 ,\dots) \not\in G_p$. Then there exists the minimal $i>1$ such that $z_i \not= z_1^{\gamma (i)}$. Set $H_i =\{ (z, z^{\gamma (2)},\dots,z^{\gamma (i)} ),\, z\in \mathbb{T} \}$. Then $H_i$ is closed in $\mathbb{T}^i$ and $\omega_0' =(z_1,\dots,z_i) \not\in H_i$. Let $\pi_i$ be the natural projection from $\mathbb{T}_p^H$ to $\mathbb{T}^i$. It is clear that $\pi_i (G_p)\subset H_i$. If $\mathbf{n}' =(n_1,\dots,n_i) \in H_i^{\perp}$ is such that $(\mathbf{n}' , \omega_0' )\not= 1$, then $\mathbf{n} =(n_1,\dots,n_i, 0,\dots) \in G_p^{\perp}$ and $(\mathbf{n} , \omega_0 )= (\mathbf{n}' , \omega_0' )\not= 1$. Hence $G_p$ is dually closed.

{\bf 2.} {\it Let us prove that $G_p, 1<p<\infty,$ is dually embedded in $\mathbb{T}_p^H$ and $G_p^{\wedge}$ is algebraically isomorphic to $\mathbb{Z}$.}

As it was proved in \cite{AaN}, $\chi =(\omega_k )\in G_p^{\wedge}$  iff either $\omega_k =0$ for all large $k$ or  $\omega_k =a_k -1$ for all large $k$. Hence we can identify $\chi = (\omega_1, \dots, \omega_m, 0,\dots)\in G_p^{\wedge}$ with $n\in \mathbb{Z}$ in the following way (10.3, \cite{HR}): if $n=\omega_1 + \omega_2 \gamma (2) +\dots +\omega_m \gamma (m) >0$, then
\[
n\mapsto \chi = (\omega_1, \dots, \omega_m, 0,\dots) \mbox{ and}
\]
\[
-n \mapsto \chi = (a_1 -\omega_1, a_2 - \omega_2 -1,\dots, a_m -\omega_m -1, a_{m+1} -1, a_{m+2} -1,\dots).
\]
Therefore $G_p^{\wedge} =\mathbb{Z} =\mathbb{T}^{\wedge}$ and
\[
(n,z) = z^n , \; \forall n\in\mathbb{Z}, \, z\in G_p .
\]
Hence we can extend every $n\in G_p^{\wedge}$ to a character of $(\mathbb{T}_p^H )^{\wedge}$, for example, in the following way
\[
n\mapsto \mathbf{n} = (n,0,\dots) \mbox{ and } (\mathbf{n}, \omega ) =z_1^n , \mbox{ where } \omega = (z_1 , z_2 ,\dots).
\]
Hence $G_p$ is a dually embedded subgroup of $\mathbb{T}_p^H$.

{\bf 3.} {\it Let us compute $G_p^{\perp}$}.
By definition, we have
\[
G_p^{\perp} =\{ \mathbf{n} =(n_1, n_2,\dots , n_s, 0,\dots) : \; z^{\sum_{k=1}^s n_k\gamma (k)} =1, \; \forall z\in G_p \}.
\]
Since $Q$ is dense in $G_p$ and $\mathbb{T}$, we have
\[
\mathbf{n} \in G_p^{\perp} \mbox{ if and only if } \sum_{k=1}^s n_k\gamma (k) =0.
\]

{\bf 4.} {\it Let us prove that $(\mathbb{T}_p^H)^{\wedge} /G_p^{\perp} =
(\mathbb{Z}, {\boldsymbol \gamma})$.}

Denote by $\pi_0 : (\mathbb{T}_p^H)^{\wedge} \mapsto (\mathbb{T}_p^H)^{\wedge} /G_p^{\perp}$ the natural homomorphism.
By item 3, if $\chi =(n_1, n_2,\dots , n_s, 0,\dots)\in (\mathbb{T}_p^H)^{\wedge}$, then $\chi + G_p^{\perp} =\chi' +G_p^{\perp}$, where $\chi' =(n,0,\dots)$ with $n=\sum_{k=1}^s n_k\gamma (k)$. Thus
$\pi_0 (e_s) = \gamma_s$. Hence, by Lemma 3 \cite{Ga3}, $(\mathbb{T}_p^H)^{\wedge} /G_p^{\perp} =
(\mathbb{Z}, {\boldsymbol \gamma})$.

In what follows we need some notations. Set
\[
A(k,m) =\left\{ n_1 \gamma (r_1) +\dots +n_s \gamma (r_s) | m\leq r_1 <\dots < r_s ,  \sum_{i=1}^s | n_i | \leq k \right\},
\]
\[
A^H (k,m) =\left\{ n_1 e_{r_1} +\dots +n_s e_{r_s} | m\leq r_1 <\dots < r_s , n_i \in \mathbb{Z}, \sum_{i=1}^s | n_i | \leq k \right\}.
\]
Then $A(k,m)$ is a subset of $(\mathbb{T}_p^H)^{\wedge} /G_p^{\perp}$, and $A^H (k,m)$ is a subset of $ \left( \mathbb{T}_0^H \right)^{\wedge}$.

Let us consider the embedding $S_p : G_p \to \mathbb{T}_p^H$. By item 2, $G_p$ is dually embedded. Thus $S_p^{\ast} : (\mathbb{T}_p^H )^{\wedge} \to G_p^{\wedge}$ is surjective. So $S_p^{\ast\ast} :G_p^{\wedge\wedge} \to (\mathbb{T}_p^H )^{\wedge\wedge} =\mathbb{T}_0^H$ is a continuous monomorphism and $\phi_p : (\mathbb{T}_p^H)^{\wedge} /G_p^{\perp} \mapsto G_p^{\wedge}$ is a continuous isomorphism.

{\bf 5.} {\it Let us prove that for every $\delta >0$ there exists $k>0$ such that $U_{\delta}^{\triangleright} \subset \phi_p (A(k,0))$.}

Let $m \in U_{\delta}^{\triangleright}$. If $m>0$, then there exists the following unique decomposition of $m$
\begin{equation} \label{20}
m= \omega_1 \gamma (1) +\dots + \omega_s \gamma (s), \mbox{ where } 0\leq\omega_k < a_k.
\end{equation}

a) {\it Let us prove that there exist $C_1 =C_1 (\delta) >0$ such that for all $m\in U_{\delta}^{\triangleright}$ with decomposition {\rm (\ref{20})} we have $\omega_k < C_1$.}

Indeed, assume the converse and there exist $m_k \in U_{\delta}^{\triangleright}, m_k = \sum_n \omega_n^k \gamma (n)$, and an index $n_k$ such that $\omega_{n_k}^k \to \infty$. We will assume that $\omega_{n_k}^k >2$ and identify $\mathbb{T}$ with $[0;1)$. Set
\[
x_k = \frac{\varepsilon_{n_k}}{a_1 a_2\dots a_{n_k}} , \mbox{ where } \varepsilon_{n_k} = \left[ \frac{3 a_{n_k}}{2 \omega_{n_k}^k} \right].
\]
Then $\gamma (n) x_k ({\rm mod 1})= \frac{\varepsilon_{n_k}}{a_n a_{n+1}\dots a_{n_k}}$ if $n\leq n_k$, and
$\gamma (n) x_k ({\rm mod 1}) =0$ if $n>n_k$.
Therefore, by (\ref{2}), we have
\[
\| x_k \|^p_p \leq \sum_{n=1}^{n_k} \frac{2^p \pi^p \varepsilon_{n_k}^p}{(a_n a_{n+1}\dots a_{n_k})^p} <\frac{2^{p+1} \pi^p \varepsilon_{n_k}^p}{a_{n_k}^p} < 2\cdot 3^p \pi^p \frac{1}{(\omega^k_{n_k})^p} <\delta^p
\]
if $(\omega^k_{n_k})^p > \frac{2\cdot 3^p \pi^p }{\delta^p}$. Hence for some $k_0$, we have $x_k \in U_{\delta} , \forall k\geq k_0$. Since
\[
\frac{1}{2\pi i} {\rm Arg} (m_k , x_k) ({\rm mod 1})= \sum_{n=1}^{n_k} \omega_n^k \frac{\varepsilon_{n_k}}{a_n a_{n+1}\dots a_{n_k}} =
\]
\[
\frac{\varepsilon_{n_k}}{a_{n_k}} \left( \omega_{n_k}^k + \sum_{n=1}^{n_k -1} \frac{\omega_n^k }{a_n \dots a_{n_k -1}} \right) ,
\]
then
\[
\frac{\omega_{n_k}^k \varepsilon_{n_k}}{a_{n_k}} <\frac{1}{2\pi i} {\rm Arg} (m_k , x_k) ({\rm mod 1}) <\frac{(\omega_{n_k}^k +2)\varepsilon_{n_k}}{a_{n_k}} \mbox{ and }
\]
\[
\frac{1}{2\pi i} {\rm Arg} (m_k , x_k) ({\rm mod 1}) \to \frac{3}{2} ({\rm mod 1}) = \frac{1}{2}.
\]
Hence $m_k \not\in U_{\delta}^{\triangleright}$. It is a contradiction.

b) {\it Let us prove that there exists an integer $M=M(\delta)$ such that at most $M$ coefficients $\omega_k$ are not equal to 0.}

Denote by $M(n)$ the number of all nonzero coefficients $\omega_i$ of $n$ in the decomposition (\ref{20}). Denote by $Z^A (n)$ the set of all coefficients $\omega_i$ of $n$ which are equal to $A$. Set $M^A (n)$ is the cardinality of $Z^A (n)$. By a), we need to prove that $M^A (n)$ is bounded on $U_{\delta}^{\triangleright}$ for every $A$.

Let us assume the converse and there exists a sequence $\{ m_k \} \subset U_{\delta}^{\triangleright}$ such that $M^A (m_k )\to\infty$. Choose $k'$ such that (we remaind that $a_k\to \infty$)
\begin{equation} \label{22}
2^{2p+1} \pi^p C_1 \sum_{k=k'}^{\infty} \frac{1}{k^p} <\delta^p \mbox{ and } a_k >100, \, \forall k\geq k' .
\end{equation}
Choose a sequence $2^p C_1 k' < T_0 < T_1 \dots$ such that
\begin{equation} \label{23}
\frac{2}{5} < \sum_{l= T_i +1}^{T_{i +1}} \frac{A}{l} <\frac{3}{5} , \quad \forall i=0,1,\dots
\end{equation}
and choose $k_0 >10 k'$ such that $a_k > 10T_1$ for all $k\geq k_0$.

Choose $m_1\in U_{\delta}^{\triangleright}$ such that $M^A (m_1)> T_1 + k_0$ and set $k'_1 >k_0$ is the maximal index of nonzero $\omega_i$ in the decomposition (\ref{20}) of $m_1$. Let us denote by $k_0 < l_1 <\dots <l_{T_1}$ the first $T_1$ indexes such that $\omega_{l_i} = A$. Choose $k_1 >k'_1$ such that $a_k >10^2 T_2$ for all $k\geq k_1$.

Choose $m_2\in U_{\delta}^{\triangleright}, m_2 \not= m_1,$ such that $M^A (m_2)> T_2 + k_1$ and set $k'_2 >k_1$ is the maximal index of nonzero $\omega_i$ in the decomposition (\ref{20}) of $m_2$. Let us denote by $k_1 < l_{T_1 +1} <\dots <l_{T_2}$ the first $T_2 -T_1$ indexes such that $\omega_{l_i} = A$. Choose $k_2 >k'_2$ such that $a_k >10^3 T_3$ for all $k\geq k_2$. And so on. Remark that, by our choosing of $l_k$,
\begin{equation} \label{24}
a_{l_{T_{k-1} +n}} > 10^k T_k , \quad k=1,2\dots, \; 1\leq n\leq T_k - T_{k-1} .
\end{equation}
Set
\[
x_k = \sum_{n=1}^{T_k -T_{k-1}} \left[ \frac{a_{l_{T_{k-1} +n}}}{T_{k-1} +n}\right] \cdot \frac{1}{a_1 a_2 \dots a_{l_{T_{k-1} +n}}} , \quad k=1,2\dots
\]
Then we have: for $1\leq s \leq l_{T_{k-1} +1}$
\[
\gamma (s) x_k ({\rm mod 1}) =
\]
\[
\frac{1}{a_{s} a_{s+1} \dots a_{l_{T_{k-1} +1} -1}} \sum_{n=1}^{T_k -T_{k-1}} \left[ \frac{a_{l_{T_{k-1} +n}}}{T_{k-1} +n}\right] \cdot \frac{1}{a_{l_{T_{k-1} +1}} \dots a_{l_{T_{k-1} +n}}} =
\]
\begin{equation} \label{25}
\frac{1}{a_{s} a_{s+1} \dots a_{l_{T_{k-1} +1} -1}} \left( \left[ \frac{a_{l_{T_{k-1} +1}}}{T_{k-1} +1}\right] + \theta_s^k \right) \cdot \frac{1}{a_{l_{T_{k-1} +1}} },
\end{equation}
where $0<\theta_s^k <1$; for $l_{T_{k-1} +r} < s \leq l_{T_{k-1} +r+1}, 0<r<T_k -T_{k-1},$
\[
\gamma (s) x_k ({\rm mod 1}) =
\]
\[
\frac{1}{a_{s} a_{s+1} \dots a_{l_{T_{k-1} +r+1} -1}} \sum_{n=r+1}^{T_k -T_{k-1}} \left[ \frac{a_{l_{T_{k-1} +n}}}{T_{k-1} +n}\right] \cdot \frac{1}{a_{l_{T_{k-1} +r+1}} \dots a_{l_{T_{k-1} +n}}} =
\]
\begin{equation} \label{26}
\frac{1}{a_{s} a_{s+1} \dots a_{l_{T_{k-1} +r+1} -1}} \left( \left[ \frac{a_{l_{T_{k-1} +r+1}}}{T_{k-1} +r+1}\right] + \theta_s^k \right) \cdot \frac{1}{a_{l_{T_{k-1} +r+1}} },
\end{equation}
where $0<\theta_s^k <1$; and $\gamma (s) x_k ({\rm mod 1}) =0$ for $l_{T_{k}} < s $.

Since $a_k \geq 2$ and $p>1$, for every $q\geq 1$, we have
\begin{equation} \label{27}
\sum_{s=1}^q \frac{1}{(a_s a_{s+1}\dots a_q)^p} \leq \sum_{s=1}^q \frac{1}{a_s a_{s+1}\dots a_q} <2.
\end{equation}

Since $a_k \geq 2, p>1$ and $T_0 >k'$, by (\ref{2}) and (\ref{22})-(\ref{27}), we have
\[
\| x_k \|_p^p < 2^p \pi^p\sum_{n=1}^{T_k -T_{k-1}} 2\left( 2\left[ \frac{a_{l_{T_{k-1} +n}}}{T_{k-1} +n}\right] \right)^p \cdot \frac{1}{ a^p_{l_{T_{k-1} +n}}} <
\]
\[
2^{2p+1} \pi^p \sum_{n=1}^{T_k -T_{k-1}} \left( \frac{1}{T_{k-1} +n} \right)^p <\delta^p .
\]
Hence $x_k \in U_{\delta}$ for all $k$.

For $m_k =\sum_n \omega_n^k \gamma (n)$ we denote by
\begin{itemize}
\item[-] $B^k_0$ is the set of indexes $1<s< l_{T_{k-1}}$ such that $\omega^k_s \not= 0$.
\item[-] $B^k_r, 0<r<T_k -T_{k-1},$ is the set of indexes $l_{T_{k-1} +r}<s< l_{T_{k-1}+r+1}$ such that $\omega^k_s \not= 0$.
\item[-] $B^k_{T_k -T_{k-1}}$ is the set of indexes $l_{T_{k}}<s$ such that $\omega^k_s \not= 0$.
\end{itemize}
Then we can represent $m_k$ in the form
\[
m_k = \sum_{n=1}^{T_k -T_{k-1}} A \gamma \left(l_{T_{k-1} +n}\right) +
\sum_{r=0}^{T_k -T_{k-1}} \sum_{s\in B^k_r}  \omega_s^k \gamma (s).
\]
Then, by (\ref{25}) and (\ref{26}), we have ($0<\theta_s^k <1$)
\[
A \gamma \left(l_{T_{k-1} +n}\right) x_k ({\rm mod} 1) =
\]
\begin{equation} \label{28}
\frac{A}{T_{k-1} +n} +
\frac{A}{a_{l_{T_{k-1} +n}}} \left( \left[ \frac{a_{l_{T_{k-1} +n}}}{T_{k-1} +n}\right]  -
\frac{a_{l_{T_{k-1} +n}}}{T_{k-1} +n} + \theta_s^k \right).
\end{equation}
for every $1\leq n\leq T_k -T_{k-1}$. For every $0\leq r< T_k -T_{k-1}$, by (\ref{25})-(\ref{27}), we have
\[
\sum_{s\in B^k_r}  \omega_s^k \gamma (s) x_k ({\rm mod} 1) <
\]
\begin{equation} \label{29}
\sum_{s=1}^{l_{T_{k-1} +r+1} -1}
\frac{\omega_s^k}{a_{s} a_{s+1} \dots a_{l_{T_{k-1} +r+1} -2}} \cdot \frac{1}{a_{l_{T_{k-1} +r+1} -1}} <\frac{2C_1}{a_{l_{T_{k-1} +r+1} -1}},
\end{equation}
and
\begin{equation} \label{30}
\sum_{s\in B^k_{T_k -T_{k-1}}}  \omega_s^k \gamma (s) x_k ({\rm mod} 1) =0
\end{equation}

Thus, by (\ref{28})-(\ref{30}), we have $({\rm mod 1})$
\[
\left| \frac{1}{2\pi i} {\rm Arg} (m_k , x_k) - \sum_{n=1}^{T_k -T_{k-1}} \frac{A}{T_{k-1} +n}\right| <
\]
\[
\sum_{n=1}^{T_k -T_{k-1}} \frac{A}{a_{l_{T_{k-1} +n}}} +\sum_{r=0}^{T_k -T_{k-1}} \frac{2C_1}{a_{l_{T_{k-1} +r+1} -1}} <\frac{4C_1 (T_k -T_{k-1})}{10^k T_k} \to 0.
\]
Hence, by (\ref{24}), $m_k\not\in U_{\delta}^{\triangleright}$. It is a contradiction.

Since $U_{\delta}^{\triangleright}$ is symmetric, we proved the following:
for each $\delta >0$ there exists a constant $C=C(\delta)$ such that if $n\in U_{\delta}^{\triangleright}$ and
\[
|n|=\omega_1 +\omega_2 \gamma (2) +\dots+\omega_m \gamma (m) , \mbox{ then } \omega_1 +\dots +\omega_m < C.
\]

Hence $U_{\delta}^{\triangleright} \subset\phi_p (A(k,0))$ for some integer $k>0$.

{\bf 6.} {\it Let us prove that $\phi_p : (\mathbb{T}_p^H)^{\wedge} /G_p^{\perp} \mapsto G_p^{\wedge}$ is a topological isomorphism.}

Let $\delta >0$. Since $A^H (k,0)$ is compact in $\left( \mathbb{T}_0^H \right)^{\wedge}$ (see proof of Theorem \ref{t1}) and $A(k,0) = \pi_0 ( A^H (k,0) )$, then $A(k,0)$ is compact.
Therefore $\phi_p$ is a homeomorphism on $A(k,0)\supset \phi_p^{-1} (U_{\delta}^{\triangleright})$. Since $G_p^{\wedge}$ is hemicompact and any compact set in $G_p^{\wedge}$ is contained in some $U_{\delta}^{\triangleright}$, $\phi_p$ is a topological isomorphism.

{\bf 7.} {\it Let us prove that $G_p^{\wedge\wedge} =G_0$.}

 By items 4 and 6, we have $G_p^{\wedge} =  (\mathbb{T}_p^H )^{\wedge} / G_p^{\perp}=(\mathbb{Z}, {\boldsymbol \gamma})$. Thus the assertion follows from Theorem 3 \cite{Ga3}. In particular,
$G_p^{\wedge}$ does not depend on $p$.

{\bf 8.} Let $\alpha_p : G_p \mapsto G_p^{\wedge\wedge}= G_0$ be the canonical homomorphism.
Since $Q$ is dense in $G_p$ \cite{AaN} and $G_0$ (Lemma \ref{l2}), then $\alpha_p (G_p)$ is dense in $G_0$. Thus, by Proposition \ref{p3}.1, $G_p^{\wedge}$, $G_0^{\wedge}$ and $G_0$ are reflexive. In particular, $G_p^{\wedge} = G_0^{\wedge}$.  By Proposition \ref{p3}.2, $G_p$ is not locally quasi-convex. $\Box$

Note that the groups $G_p, 1<p<\infty,$ form a continual chain (under inclusion) of {\it Polish  non} locally quasi-convex totally disconnected groups with the same {\it countable reflexive} dual.

\end{document}